\documentclass[10pt,a4paper]{article}

\usepackage[hyperfootnotes=false,pdfpage layout=SinglePage]{hyperref}
\usepackage{mathptmx}
\usepackage[scaled=.92]{helvet}
\usepackage{courier}
\usepackage{graphicx}
\usepackage{marvosym}
\usepackage{pstricks,pst-eucl,pst-plot,pstricks-add,pst-text}
\usepackage{amssymb}
\usepackage{amsmath}
\usepackage{amsthm}
\usepackage{latexsym}
\usepackage{amsfonts}
\usepackage{mathrsfs}
\usepackage{setspace}

\newtheorem{theorem}{Theorem}[section]

\newtheorem{lemma}[theorem]{Lemma}

\theoremstyle{definition}
\newtheorem{definition}[theorem]{Definition}

\numberwithin{equation}{section}
\renewcommand{\emptyset}{\varnothing}

\newcommand{\R}{\ensuremath{\mathbb R}}    
\newcommand{\C}{\ensuremath{\mathbb C}}    
\newcommand{\N}{\ensuremath{\mathbb N}}    

\newcommand{\gperp}{{[\perp]}}

\newcommand{\product}{[\cdot\,,\cdot]}
\newcommand{\hproduct}{(\cdot\,,\cdot)}

\newcommand{\calD}{\mathcal D}

\newcommand{\calH}{\mathcal H}

\newcommand{\calK}{\mathcal K}         \newcommand{\frakK}{\mathfrak K}
         
\newcommand{\calM}{\mathcal M}         
\newcommand{\calN}{\mathcal N}

         \newcommand{\frakS}{\mathfrak S}
         
\newcommand{\calU}{\mathcal U}

\newcommand{\la}{\lambda}
\newcommand{\veps}{\varepsilon}
\newcommand{\vphi}{\varphi}

\renewcommand{\Im}{\operatorname{Im}}

\renewcommand{\ker}{\operatorname{ker}}
\newcommand{\ran}{\operatorname{ran}}

\newcommand{\sess}{\sigma_{\rm ess}}
\renewcommand{\sp}{\sigma_{+}}
\newcommand{\sm}{\sigma_{-}}

\newcommand{\Lra}{\Longrightarrow}

\newcommand{\upto}{\uparrow}
\newcommand{\downto}{\downarrow}

\newcommand{\ol}{\overline}
\newcommand{\ds}{\dotplus}
\newcommand{\wt}{\widetilde}

\newcommand{\dist}{\operatorname{dist}}
\newcommand{\tr}{\operatorname{tr}}

\begin{document}
\vspace*{-.3cm}
\begin{center}
\begin{spacing}{1.7}
{\LARGE\bf Variation of discrete spectra of non-negative operators in Krein spaces}
\end{spacing}

\vspace{1cm}
{\Large Jussi Behrndt, Leslie Leben and Friedrich Philipp}
\end{center}

\vspace{.5cm}
\begin{abstract}\noindent
We study the variation of the discrete spectrum of a bounded non-negative operator in a Krein space 
under a non-negative Schatten class perturbation of order $p$. It turns out that there exist so-called extended enumerations 
of discrete eigenvalues of the unperturbed and the perturbed operator, respectively, whose difference is an $\ell^p$-sequence. 
This result is a Krein space version of a theorem by T.~Kato for bounded selfadjoint operators in Hilbert spaces.
\end{abstract}

\vspace{.6cm}


\section{Introduction}
In this note we prove a Krein space version of a result by T.~Kato from \cite{k87} on the variation of the discrete spectra of bounded selfadjoint 
operators in Hilbert spaces under additive perturbations from the Schatten-von Neumann ideals $\frakS_p$. Although perturbation theory for
selfadjoint operators in Krein spaces is a well developed field, and compact, finite rank, as well as bounded 
perturbations have been studied extensively, 
only very few results exist that take into account a particular $\frakS_p$-character of perturbations.
To give an impression of the variety of perturbation results for various classes of selfadjoint operators in Krein spaces 
we refer the reader to
\cite{ajt05,bj,J82,J88,J90,J91,lmm} for compact perturbations, to \cite{ABPT08,ABT08,b07,JonasLanger79,JL83} for finite rank perturbations, 
and to \cite{AMS09,AMT09,abjt,ajt09,Jonas98,KMM07,LN83,T09} for 
(relatively) bounded and small perturbations.

Here we consider a bounded operator $A$ in a Krein space $(\calK,\product)$ which is assumed to be non-negative with 
respect to the indefinite inner product $\product$, and an additive perturbation $C$ which is also non-negative and belongs 
to some Schatten-von Neumann 
ideal $\frakS_p$, that is, 
$C$ is compact and its singular values form a sequence in $\ell^p$. 
Recall that the spectrum of a bounded non-negative operator in $(\calK,\product)$ is real.
We also assume that $0$ is not a
singular critical point of the perturbation $C$, which is a typical assumption in perturbation theory for selfadjoint operators in Krein spaces;
cf. Section 2 for a precise definition. We note that by this assumption $C$ is similar to a selfadjoint operator in a Hilbert space.
Clearly, the non-negativity and compactness of $C$ imply that the bounded operator 
$$
B := A + C
$$
is also non-negative in $(\calK,\product)$ and its essential spectrum coincides with that of $A$, whereas
the discrete eigenvalues of $A$ and their multiplicity are in general not stable under the perturbation $C$. 
Hence, it is particularly interesting to prove qualitative and quantitative results on the discrete spectrum.
Our main objective here is to 
compare the discrete spectra of $A$ and $B$. For that we make use of
the following notion from \cite{k87}: Let $\Delta\subset\R$ be a finite union of open intervals. A sequence  $(\alpha_n)$ is said to be an
{\it extended enumeration of discrete eigenvalues of $A$ in $\Delta$} if every discrete eigenvalue of $A$ in $\Delta$ 
with multiplicity $m$ appears exactly $m$-times
in the values of $(\alpha_n)$ and all other values $\alpha_n$ are boundary points of the essential spectrum of $A$ in $\ol\Delta\subset\R$.
An extended enumeration of discrete eigenvalues of $B$ in $\Delta$ is defined analogously.
The following theorem is the main result of this note.

\begin{theorem}\label{t:main}
Let $A$ and $B$ be bounded non-negative operators in a Krein space $(\calK,\product)$ such that $B=A+C$, where $C\in\frakS_p(\calK)$ 
is non-negative, $0$ is not a singular critical point of $C$ and $\ker C = \ker C^2$. Then for each finite union of open intervals $\Delta$ with $0\notin\ol\Delta$ there exist extended enumerations $(\alpha_n)$ and $(\beta_n)$ of the discrete eigenvalues of $A$ and $B$ in $\Delta$, respectively, such that
$$
(\beta_n - \alpha_n)\,\in\,\ell^p.
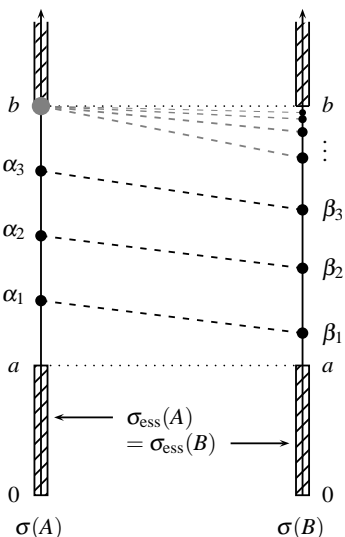
$$
\end{theorem}

\begin{minipage}{0.42\textwidth}

\scalebox{0.86}{
\mbox{
\begin{pspicture}(3,-1.3)(7.8,7.2)

\pspolygon[linecolor=white,fillstyle=hlines](2.9,5.5)(2.9,6.8)(3.1,6.8)(3.1,5.5) 
\pspolygon[linecolor=white,fillstyle=hlines](6.9,5.5)(6.9,6.8)(7.1,6.8)(7.1,5.5) 
\psline{-}(2.9,5.5)(2.9,6.8)
\psline{-}(3.1,6.8)(3.1,5.5)
\psline{-}(3.1,5.5)(2.9,5.5)
\psline{-}(6.9,6.8)(6.9,5.5)
\psline{-}(6.9,5.5)(7.1,5.5)
\psline{-}(7.1,5.5)(7.1,6.8)
\pspolygon[linecolor=black,fillstyle=hlines](2.9,-0.5)(2.9,1.5)(3.1,1.5)(3.1,-0.5) 
\pspolygon[linecolor=black,fillstyle=hlines](6.9,-0.5)(6.9,1.5)(7.1,1.5)(7.1,-0.5) 

\psline{->}(3,-0.5)(3,7) 
\psline{->}(7,-0.5)(7,7) 

\psline[linestyle=dotted](3,1.5)(7,1.5) \put(2.5,1.4){$a$} \put(7.3,1.4){$a$}
\psline[linestyle=dotted](3,5.5)(7,5.5) \put(2.5,5.4){$b$} \put(7.3,5.4){$b$}
\put(2.5,-0.6){$0$} \put(7.3,-0.6){$0$}

\qdisk(3,2.5){2.5pt}
\qdisk(3,3.5){2.5pt}
\qdisk(3,4.5){2.5pt}
\psset{linecolor=gray}  
\qdisk(3,5.5){4pt}    
\psset{linecolor=black} 

\put(2.4,2.5){$\alpha_1$}
\put(2.4,3.5){$\alpha_2$}
\put(2.4,4.5){$\alpha_3$}

\psline[linestyle=dashed,dash=2.8pt 3.5pt](3,2.5)(7,2)
\psline[linestyle=dashed,dash=2.8pt 3.5pt](3,3.5)(7,3)
\psline[linestyle=dashed,dash=2.8pt 3.5pt](3,4.5)(7,3.9)
\psline[linestyle=dashed,dash=2.8pt 3.5pt,linecolor=gray](3,5.5)(7,4.7)
\psline[linestyle=dashed,dash=2.8pt 3.5pt,linewidth=0.8pt,linecolor=gray](3,5.5)(7,5.1)
\psline[linestyle=dashed,dash=2.8pt 3.5pt,linewidth=0.6pt,linecolor=gray](3,5.5)(7,5.3)
\psline[linestyle=dashed,dash=2.8pt 3.5pt,linewidth=0.4pt,linecolor=gray](3,5.5)(7,5.4)

\qdisk(7,2){2.5pt}
\qdisk(7,3){2.5pt}
\qdisk(7,3.9){2.5pt}
\qdisk(7,4.7){2.5pt}
\qdisk(7,5.1){2.2pt}
\qdisk(7,5.3){1.8pt}
\qdisk(7,5.4){1.5pt} 

\put(7.3,1.9){$\beta_1$}
\put(7.3,2.9){$\beta_2$}
\put(7.3,3.8){$\beta_3$}
\put(7.3,4.65){$\vdots$}

\put(2.6,-1.1){$\sigma(A)$}
\put(6.6,-1.1){$\sigma(B)$}
\put(4.3,0.6){$\sess(A)$}
\psline{->}(4.1,0.7)(3.2,0.7)
\put(4.3,0.2){$= \sess(B)$}
\psline{->}(5.9,0.3)(6.8,0.3)

\end{pspicture}
}}
\end{minipage}
\begin{minipage}{0.5\textwidth}
The adjacent figure illustrates the role of extended enumerations in Theorem~\ref{t:main}: We consider a gap $(a,b)\subset\R$ 
in the essential spectrum and compare the discrete spectra of $A$ and $B$ therein. Here the discrete spectrum of the 
unpertur\-bed operator $A$ in $(a,b)$ consists of the (simple) eigenvalues $\alpha_1,\alpha_2,\alpha_3$, and the 
eigenvalues $\beta_n$, $n=1,2,\dots$, of the perturbed operator $B$ accumulate to the boundary point $b\in\partial\sess(A)$. 
Therefore, in the situation of Theorem~\ref{t:main} the value $b$ is contained (infinitely many times) in the extended 
enumeration $(\alpha_n)$ of the discrete eigenvalues of $A$ in $(a,b)$.
\end{minipage}


For bounded selfadjoint operators $A$ and $B$ in a Hilbert space and an $\frakS_p$-perturbation $C$ Theorem~\ref{t:main} was proved
by T. Kato in \cite{k87}. The original proof is based on methods from analytic perturbation theory, in particular, on the properties 
of a family of real-analytic functions describing the discrete eigenvalues and eigenprojections of the operators $A(t)=A+tC$, $t\in\R$.
Our proof follows the lines of Kato's proof, but in the Krein space situation some nontrivial additional arguments 
and adaptions are necessary. In particular,
we apply methods from \cite{lmm} to show that the non-negativity assumptions on $A$ and $C$ yield 
uniform boundedness of the spectral projections of $A(t)$, $t\in[0,1]$, corresponding to positive and negative intervals,
respectively. 
The non-negativity assumptions on $A$ and $C$ also enter in the construction and properties of the real-analytic functions associated with
the discrete eigenvalues of $A(t)$. 

Besides the introduction this note consists of two further sections. In Section~2 we recall some definitions and 
spectral properties of non-negative operators in Krein spaces. Section~3 contains the proof of our main result Theorem~\ref{t:main}.
As a preparation, we discuss the properties of the family of real-analytic functions describing the
eigenvalues and eigenspaces of $A(t)$ in Lemma 3.1 and show a result on the uniform definiteness 
of certain spectral subspaces of $A(t)$ in Lemma 3.2. Finally, by modifying and following some of the arguments and estimates in \cite{k87} we complete 
the proof of our main result.

\section{Preliminaries on non-negative operators in Krein spaces}
Throughout this paper let $(\calK,\product)$ be a Krein space. For a detailed study of Krein spaces and operators therein we 
refer to the monographs \cite{ai} and \cite{bo}. For the rest of this section let $\|\cdot\|$ be a Banach space norm with 
respect to which the inner product $\product$ is continuous. All such norms are equivalent, see \cite{ai}. For closed 
subspaces $\calM$ and $\calN$ of $\calK$ we denote by $L(\calM,\calN)$ the set of all bounded and everywhere defined 
linear operators from $\calM$ to $\calN$. As usual, we write $L(\calM) := L(\calM,\calM)$. 

Let $T\in L(\calK)$.
The adjoint of $T$, denoted by $T^+$, is defined by 
$$
[Tx,y] = [x,T^+y]\quad\text{for all }x,y\in\calK.
$$
The operator $T$ is called {\it selfadjoint} in $(\calK,\product)$ (or $\product$-{\it selfadjoint}) 
if $T = T^+$. 
Equivalently, $[Tx,x]\in\R$ for all $x\in\calK$. We mention that the spectrum of a selfadjoint operator in a Krein space 
is symmetric with respect to the real axis but in general not contained in $\R$. 

The following definition of spectral points of positive and negative type is from \cite{lmm}.

\begin{definition}
Let $A\in L(\calK)$ be a selfadjoint operator. A point $\la\in\sigma(A)\cap\R$ is called 
a {\em spectral point of positive type {\rm (}negative type{\rm )} of} $A$ if for each sequence 
$(x_n)\subset\calK$ with $\|x_n\| = 1$, $n\in\N$, and $(A - \la)x_n\to 0$ as $n\to\infty$ we have
$$
\liminf_{n\to\infty}\,[x_n,x_n] > 0\qquad\Big(\limsup_{n\to\infty}\,[x_n,x_n] < 0,\text{ respectively}\Big).
$$
The set of all spectral points of positive {\rm (}negative{\rm )} type of 
$A$ is denoted by $\sp(A)$ {\rm (}$\sm(A)$, respectively{\rm )}. A set $\Delta\subset\R$ is 
said to be of {\em positive type {\rm (}negative type{\rm )}} with respect to $A$ if each spectral 
point of $A$ in $\Delta$ is of positive type {\rm (}negative type, respectively{\rm )}.
\end{definition}

A closed subspace $\calM\subset\calK$ is called {\it uniformly positive} ({\it uniformly negative}) 
if there exists $\delta > 0$ such that $[x,x]\ge\delta\|x\|^2$ ($[x,x]\le -\delta\|x\|^2$, respectively)
holds for all $x\in\calM$. Equivalently, $(\calM,\product)$ ($(\calM,-\product)$, respectively) is a Hilbert 
space. For a bounded selfadjoint operator $A$ in $\calK$ it follows directly from the definition of $\sp(A)$ and $\sm(A)$ that 
an isolated eigenvalue $\la_0\in\R$ of $A$ is of positive type (negative type) 
if and only if $\ker(A - \la_0)$ is uniformly positive (uniformly negative, respectively).

A selfadjoint operator 
$A\in L(\calK)$ is called {\it non-negative} if
$$
[Ax,x]\ge 0\quad\text{for all }x\in\calK.
$$
The spectrum of a bounded non-negative operator $A$ is a compact subset of $\R$ and
\begin{equation}\label{guteformel}
\sigma(A)\cap\R^\pm\subset\sigma_\pm(A)
\end{equation}
holds, see \cite{l}. The {\it discrete spectrum} $\sigma_d(A)$ of $A$ consists of the isolated eigenvalues of $A$ with 
finite multiplicity. The remaining part of $\sigma(A)$ is the {\it essential spectrum} of the nonnegative operator $A$ 
and is denoted by $\sess(A)$. Observe that $\sess(A)$ coincides with the set of $\lambda$ such that $A-\lambda$ is not
a Semi-Fredholm operator.
Recall that the non-negative operator $A$ admits a spectral function $E$ on $\R$ with a possible singularity at zero, 
see \cite{l}. The spectral projection $E(\Delta)$ is defined for all Borel sets $\Delta\subset\R$ with $0\notin\partial\Delta$ 
and is selfadjoint. Hence,
$$
\calK = E(\Delta)\calK\,[\ds]\,(I - E(\Delta))\calK,
$$
which implies that $(E(\Delta)\calK,\product)$ is itself a Krein space. 
For $\Delta\subset\R^\pm$, $0\notin\ol\Delta$, 
the spectral subspace $(E(\Delta)\calK,\pm\product)$ is a Hilbert space; cf. \cite{l,lmm} and \eqref{guteformel}.

The point zero is called a {\it critical point} of a non-negative operator $A\in L(\calK)$ if 
$0\in\sigma(A)$ is neither of positive nor negative type. If zero is a critical point of $A$, 
it is called {\it regular} if $\|E([-\frac 1 n,\frac 1 n])\|$, $n\in\N$, is uniformly bounded, i.e.\ if 
zero is not a singularity of the spectral function $E$. Otherwise, the critical point zero is called {\it singular}. 
It should be noted that the non-negative operator $A\in L(\calK)$ is (similar to) a selfadjoint operator 
in a Hilbert space if and only if zero is not a singular critical point of $A$.

\section{Proof of Theorem 1.1}
Throughout this section let $A$, $B$ and $C$ be bounded non-negative operators in the Krein space $(\calK,\product)$ 
as in Theorem \ref{t:main}. By assumption $0$ is not a  singular critical point of $C$ and $C\in\frakS_p(\calK)$. 
In order to prove Theorem \ref{t:main} we consider the analytic operator function
$$
A(z) := A + zC,\quad z\in\C.
$$
Note that $A(t)$ is non-negative for $t\ge 0$. Moreover, since $C$ is compact, the essential spectrum of $A(z)$ 
does not depend on $z$ and hence
\begin{equation}\label{sessz}
\sess(A)=\sess(B)=\sess(A(z)),\qquad z\in\C.
\end{equation}
The following lemma describes the evolution of the {\it discrete} spectra of the operators $A(t)$, $t\ge 0$.

\begin{lemma}\label{l:evf}
Assume that $\sigma_d(A(t_0))\not=\emptyset$ for some $t_0\geq 0$. Then
there exist intervals $\Delta_j\subset \R_0^+$, $j=1,\dots, m$ or $j\in\N$, 
and real-analytic functions 
\begin{equation*}
\la_j(\cdot) : \Delta_j \to \R_0^+\quad\text{and}\quad
E_j(\cdot) : \Delta_j \to L(\calK),
\end{equation*}
such that the following holds.
\begin{enumerate}
\item[{\rm (i)}]  The sets $\Delta_j$ are $\R_0^+$-open intervals {\rm (}in $\R_0^+${\rm )} which are maximal with 
respect to {\rm (ii)--(vi)} below.
\item[{\rm (ii)}]  For each $t\ge 0$ we have 
\begin{align*}
\sigma_d\big( A(t)\big)\cap\R^+ = \bigl\{\la_j(t) : j\in\N \text{ such that } t \in \Delta_j \text{ and } \la_j(t)\ne 0\bigr\}.
\end{align*}
\item[{\rm (iii)}] For all $j$ and $t\in\Delta_j$ the set $\{k\in\N : \la_k(t) = \la_j(t)\}$ is finite and
$$
\sum_{k:\la_k(t) = \la_j(t)} E_k(t)
$$
is the $\product$-selfadjoint projection onto $\ker(A(t)-\la_j(t))$.
\item[{\rm (iv)}]  For all $j$ the value
$$
m_j := \dim E_j(t)\calK,\quad t\in\Delta_j,
$$
is constant.
\item[{\rm (v)}]   For all $j$ and $t\in\Delta_j$ there exists an orthonormal basis $\{x_i^j (t)\}_{i=1}^{m_j}$ of the 
Hilbert space $(E_j(t)\calK,\product)$, such that the functions $x_i^j(\cdot): \Delta_j\to\calK$ are real-analytic and the differential equation
\begin{equation}\label{e:DGL}
\la_j'(t) = \frac{1}{m_j}\sum_{k=1}^{m_j}\left[Cx_k^j(t), x_k^j(t)\right]\,\ge\,0
\end{equation}
holds. In particular, $\la_j'(t) = 0$ implies $E_j(t)\calK\subset\ker C$.
\item[{\rm (vi)}] Let $\R^+\setminus\sess(A) = \dot{\bigcup}_{n}\,\calU_n$ with mutually disjoint open intervals 
$\calU_n\subset\R^+$. For every $j$ there exists $n\in\N$ such that
\begin{align*}
\la_j(t)\in\calU_n\,\text{ for all }t\in\Delta_j&, \quad\text{ if } 0\notin\partial\calU_n,\\
\la_j(t)\in\calU_n\cup\{0\}\,\text{ for all }t\in\Delta_j&,\quad \text{ if } 0\in\partial\calU_n.
\end{align*}
If $\sup\Delta_j < \infty$, then $\sup\calU_n < \infty$ and $\lim\limits_{t\upto\sup\Delta_j} \la_j(t) = \sup\calU_n$. Moreover,
\begin{align*}
\lim_{t\downto\inf\Delta_j}\la_j(t) = \inf\calU_n&, \quad \text{ if $\Delta_j$ is open},\\
\lim_{t\downto 0}\la_j(t)\in\calU_n\cup\{\inf\calU_n\}&, \quad \text{ if $\Delta_j= [0,\sup \Delta_j)$}.
\end{align*}
\end{enumerate}
\end{lemma}


\hspace*{0.7cm}
\scalebox{0.9}{
\mbox{
\begin{pspicture}(-1,-2)(11,4.7)

\pspolygon[linecolor=white,fillstyle=hlines,hatchcolor=gray](-0.15,-0.4)(-0.15,1)(0.15,1)(0.15,-0.4)
\psline[linecolor=gray](-0.15,-0.4)(-0.15,1)
\psline[linecolor=gray](-0.15,1)(0.15,1)
\psline[linecolor=gray](0.15,1)(0.15,-0.4)
\pspolygon[linecolor=white,fillstyle=hlines,hatchcolor=gray](-0.15,4)(-0.15,4.4)(0.15,4.4)(0.15,4)
\psline[linecolor=gray](-0.15,4.4)(-0.15,4)
\psline[linecolor=gray](-0.15,4)(0.15,4)
\psline[linecolor=gray](0.15,4)(0.15,4.4)
\put(-1.4,0.4){$\sess(A)$}

\pscurve[linestyle=dotted,linewidth=1.5pt,dotsep=1.5pt](1.4,1)(3.5,2.3)(4.6,4)
\pscurve[](0,1.6)(4.9,2.8)(7,4)
\pscurve[linecolor=lightgray](0,1)(7,2.3)(9,2.6)
\pscurve[linestyle=dashed](1,1)(3,2.6)(7,3.3)(8.5,4)

\psline{->}(-0.5,0)(9.3,0) 
\psline{->}(0,-0.5)(0,4.6) 
\psline[linestyle=dotted](0,1)(9,1) 
\psline[linestyle=dotted](0,4)(9,4) 
\psline[linestyle=dotted](6.5,0)(6.5,4) 

\put(6.4,-0.5){$1$}
\put(9.1,-0.5){$t$}
\put(-0.5,0.9){$a$}
\put(-0.5,3.9){$b$}

\put(9.4,3.5){$\lambda_1(t)$} 
\psline(10.4,3.6)(10.8,3.6)
\put(9.4,2.8){$\lambda_2(t)$} 
\psline[linestyle=dashed](10.4,2.9)(10.8,2.9)
\put(9.4,2.1){$\lambda_3(t)$}
\psline[linestyle=dotted,linewidth=1.5pt,dotsep=1.5pt](10.4,2.2)(10.8,2.2)
\put(9.4,1.4){$\lambda_4(t)$} 
\psline[linecolor=lightgray](10.4,1.5)(10.8,1.5)

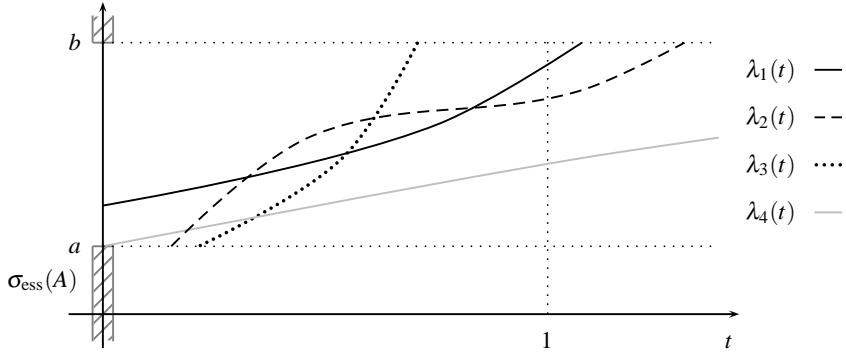
\put(-0.7,-1){Typical situation for the evolution of the discrete eigenvalues of the } 
\put(-0.7,-1.4){operator function $A(\cdot)$ in a gap $(a,b)\subset\R$ of the essential spectrum. }
\end{pspicture}
}}


\begin{proof}
The proof is based on the analytic perturbation theory of discrete eigenvalues, cf.\ \cite[Chapter\ II~and~VII]{k}, \cite{ba} and 
\cite{k87}. We fix some $t_0\ge 0$ for which an eigenvalue $\la_0\in\sigma_d(A(t_0))\cap\R^+$ exists and set 
$M(t_0) := \ker(A(t_0) - \la_0)$. Due to the non-negativity of $A$ and $C$ and since $\la_0 > 0$, the inner 
product space $(M(t_0),\product)$ is a (finite-dimensional) Hilbert space; cf. \eqref{guteformel}. Therefore, the decomposition
$$
\calK = M(t_0)\,[\ds]\,M(t_0)^\gperp
$$
reduces the operator $A(t_0)$. As in \cite[ch.\ VII]{k} one shows that for $z$ in a (complex) neighborhood $\calD$ of $t_0$ there 
exists an analytic operator function $U(\cdot) : \calD\to L(\calK)$ with $U(z)^{-1} = U(\ol z)^+$, $U(t_0) = I$ 
and such that $M(t_0)$ is $U(z)^{-1}A(z)U(z)$-invariant, $z \in \calD$. 
Hence, there exist a finite number of (possibly multivalued)  
analytic functions $\la_k(\cdot)$ describing the eigenvalues of $B(z) := U(z)^{-1}A(z)U(z)|M(t_0)$ for $z\in\calD$, see, e.g., \cite{ba}. 
Since for real $t\in\calD$ the operator $B(t)$ is selfadjoint in the Hilbert space $(M(t_0),\product)$ it follows from \cite[Theorem~II-1.10]{k} 
that the functions $\la_k(\cdot)$ are in fact single-valued. 
The same is true for the eigenprojection functions $E_k(\cdot)$,
$$
E_k(z) = -\frac 1 {2\pi i}\,\int_{\Gamma_k(z)}\,(A(z) - \la)^{-1}\,d\la,\quad z\in\calD,
$$
where $\Gamma_k(z)$ is a small circle with center $\la_k(z)$. Now a continuation argument 
implies that there exist functions $\la_j(\cdot)$, $E_j(\cdot)$
with the properties (i)--(iv) and (vi), cf. \cite{k87}.

It remains to prove (v). For this fix $j\in\N$ and $t_0\in\Delta_j$. Similarly as above there exists a function  
$U_j(\cdot): \Delta_j \to E_j(t_0)\calK$ with $U_j(t)^+ = U_j(t)^{-1}$, $U_j(t_0)=I$, and $E_j(t) = U_j(t)^+E_j(t_0)U_j(t)$ 
for every $t\in\Delta_j$. We choose an orthonormal basis $\{x_1, ..., x_{m_j}\}$ of the $m_j$-dimensional Hilbert 
space $(E_j(t_0)\calK, \product)$ and define
\begin{align*}
x_k(t) := U_j(t)x_k, \quad t \in \Delta_j,\; k=1, ..., m_j.
\end{align*}
For every $t \in \Delta_j$ the set $\{x_1(t), ..., x_{m_j}(t) \}$ is an orthonormal basis of the subspace $(E_j(t)\calK,\product)$, since for $k,l\in\{1,\ldots,m_j\}$ we have
\begin{align*}
[x_k(t), x_l(t)] = [U_j(t)x_k, U_j(t)x_l] = [x_k, x_l] = \delta_{kl}.
\end{align*}
Let $k\in\{1, ..., m_j\}$. Then
$$
[x_k'(t), x_k(t)] + [x_k(t), x_k'(t)]  = \frac{d}{dt}[x_k(t),x_k(t)] = 0
$$
and hence
\begin{align*}
\la_j'(t)
&= \frac{d}{dt} [\lambda_j(t)x_k(t), x_k(t)] = \frac{d}{dt} [A(t)x_k(t), x_k(t)]\\
&= [C x_k(t), x_k(t)] + [A(t)x_k'(t), x_k(t)] + [A(t)x_k(t), x_k'(t)] \\
&= [C x_k(t), x_k(t)] + \la_j(t)[x_k'(t), x_k(t)] + \la_j(t)[x_k(t), x_k'(t)] \\
&= [C x_k(t), x_k(t)]
\geq 0.
\end{align*}
This yields \eqref{e:DGL}.
Finally, if we have $\la_j'(t)= 0$, then $[Cx_k(t), x_k(t)] =0$ holds for $k =1, ..., m_j$. Since $C$ is non-negative, the 
Cauchy-Schwarz inequality applied to the non-negative inner product $[C\cdot,\cdot]$  yields 
\begin{align*}
\|Cx_k(t)\|^2 = \left[Cx_k(t), JCx_k(t) \right] \leq \left[Cx_k(t),x_k(t) \right]^{1/2} \left[CJCx_k(t), JCx_k(t) \right]^{1/2} = 0
\end{align*}
for every $k \in \{1, ..., m_j\}$. This shows $E_j(t)\calK \subset \ker C$. 
\end{proof}

In the proof of the following lemma we make use of methods from \cite{lmm} in order to show the uniform 
definiteness of a family of spectral subspaces of $A(t)$.

\begin{lemma}\label{l:lmm}
Let  $E_{A(t)}$ be the spectral function of the non-negative operator $A(t)$, $t\ge 0$, and let $a > 0$. Then there exists 
$\delta > 0$ such that for all $t\in [0,1]$ and all $x\in E_{A(t)}([a,\infty))\calK$ we have
\begin{equation}\label{e:lmm}
[x,x]\,\ge\,\delta\|x\|^2.
\end{equation}
\end{lemma}
\begin{proof}
Since $\max\sigma(A(t))\le b:= \Vert A\Vert +\Vert C\Vert$ for all $t\in [0,1]$, it is sufficient to
prove \eqref{e:lmm} only for $x\in E_{A(t)}([a,b])$. The proof is divided into four steps.

\vskip 0.2cm\noindent
{\bf 1.} In this step we show that there exist $\veps > 0$ and an open neighborhood $\calU$ of $[a,b]$ in $\C$ 
such that for all $t\in [0,1]$, all $\la\in\calU$ and all $x\in\calK$ we have
\begin{equation}\label{donnerwetter}
\|(A(t) - \la)x\|\le\veps\|x\|\quad\Lra\quad [x,x]\ge\veps\|x\|^2.
\end{equation}
Assume that $\veps$ and $\calU$ as above do not exist. Then there exist sequences $(t_n)\subset [0,1]$, $(\la_n)\subset\C$ and $(x_n)\subset\calK$ 
with $\|x_n\| = 1$ and $\dist(\la_n,[a,b]) < 1/n$ for all $n\in\N$, such that $\|(A(t_n) - \la_n)x_n\|\le 1/n$ and $[x_n,x_n]\le 1/n$. 
It is no restriction to assume that $\la_n\to\la_0\in [a,b]$ and $t_n\to t_0\in [0,1]$ as $n\to\infty$. Therefore,
$$
(A(t_0) - \la_0)x_n = (t_0 - t_n)Cx_n + (A(t_n) - \la_n)x_n + (\la_n - \la_0)x_n
$$
tends to zero as $n\to\infty$. But by \eqref{guteformel} we have $\la_0\in\sp(A(t_0))$ which implies $\liminf_{n\to\infty}[x_n,x_n] > 0$, 
contradicting $[x_n,x_n] < 1/n$, $n \in \N$. 

\vskip 0.2cm\noindent
{\bf 2.} In the following $\veps > 0$ and $\calU$ are fixed such that \eqref{donnerwetter} holds, and, in addition, we assume that 
$|\Im\la| < 1$ holds for all $\la\in\calU$.
Next, we verify that for all $t\in [0,1]$ 
\begin{equation}\label{e:lmm2}
\left\|(A(t) - \la)^{-1}\right\|\,\le\,\frac{\veps^{-1}}{|\Im\la|},\qquad \la\in\calU\setminus\R,
\end{equation}
holds. Indeed, for all $t\in [0,1]$, all $\la\in\calU$ and all $x\in\calK$ we either have 
$$\|(A(t) - \la)x\| > \veps\|x\|$$ 
or, by \eqref{donnerwetter},
$$
\veps|\Im\la|\|x\|^2\le |\Im\la[x,x]| = |\Im [(A(t) - \la)x,x]|\le\|(A(t) - \la)x\|\|x\|.
$$
Hence, it follows that for all $t\in [0,1]$, all $\la\in\calU$ and all $x\in\calK$ we have
$$
\|(A(t) - \la)x\|\,\ge\,\veps|\Im\la|\|x\|,
$$
which implies \eqref{e:lmm2}.

\vskip 0.2cm\noindent
{\bf 3.} In the remainder of this proof we set
$$
d := \dist([a,b],\partial\calU)\quad\text{and}\quad\tau_0 := \min\left\{\veps^2,\frac d 2\right\}.
$$ 
Let $\Delta\subset[a,b]$ be an interval of length $R\le\tau_0$ and let $\mu_0$ be the center of $\Delta$. We show that 
for all $t\in [0,1]$ the estimate
\begin{equation}\label{e:lmm3}
\big\|\big(A(t)|E_t(\Delta)\calK\big) - \mu_0\big\|\,\le\,\veps
\end{equation}
holds. For this let $B(t) := (A(t)|E_t(\Delta)\calK) - \mu_0$, $t\in[0,1]$, and note that
\begin{equation}\label{donnerwetter2}
\sigma(B(t))\,\subset\,\left[-\frac{R}{2} ,\frac{R}{2}\right]\,\subset\, (-R,R).
\end{equation}
As $R< d$, for every $\lambda\in\C\setminus\R$ with $\vert\lambda\vert<R$ we have 
$\mu_0+\la\in\calU\setminus\R$ and hence
$$
\big\|\big(B(t) - \la\big)^{-1}\big\| \le \big\|\big(A(t) - (\mu_0 + \la)\big)^{-1}\big\|\,\le\,\frac{\veps^{-1}}{|\Im\la|}
$$
by \eqref{e:lmm2}.
From \cite[Section 2(b)]{lmm} we now obtain $\|B(t)\|\le 2\veps^{-1}\,r(B(t))$, where $r(B(t))$ denotes the spectral radius of $B(t)$. 
Now \eqref{e:lmm3} follows from \eqref{donnerwetter2} and $R\leq\tau_0\leq\veps^2$. 

\vskip 0.2cm\noindent
{\bf 4.} We cover the interval $[a,b]$ with mutually disjoint intervals $\Delta_1,\ldots,\Delta_n$ of length $<\tau_0$. 
Let $\mu_j$ be the center of the interval $\Delta_j$, $j=1,\ldots,n$. From step 3 we obtain for all $t\in [0,1]$:
$$
\big\|\big(A(t)|E_{A(t)}(\Delta_j)\calK\big) - \mu_j\big\|\,\le\,\veps.
$$
Hence, by step 1 of the proof $[x_j,x_j]\ge\veps\|x_j\|^2$ for $x_j\in E_{A(t)}(\Delta_j)$, $j=1,\ldots,n$, and $t\in[0,1]$. But
$$
E_{A(t)}([a,b]) = E_{A(t)}(\Delta_1)\,[\ds]\,\ldots\,[\ds]\,E_{A(t)}(\Delta_n),
$$
and therefore with $x_j := E_{A(t)}(\Delta_j)x$, $j=1,\ldots,n$, we find that
$$ 
[x,x]\ge\veps\bigl(\Vert x_1\Vert^2+\dots+\Vert x_n\Vert^2\bigr)\geq\frac{\veps}{2^{n-1}}\,\Vert x_1+\dots + x_n\Vert^2=
\frac{\veps}{2^{n-1}}\,\Vert x \Vert^2
$$
holds for all $x\in E_{A(t)}([a,b])$ and $t\in [0,1]$, i.e. \eqref{e:lmm} holds 
with  $\delta := \veps/2^{n-1}$.
\end{proof}


\begin{proof}[Proof of Theorem \ref{t:main}]
It suffices to prove the theorem for the case that $\Delta$ is an open interval $(a,b)$ with $a > 0$. 
In the case $b < 0$ consider the non-negative operators $-A$, $-B$ and $-C$ in the Krein space $(\calK,-\product)$.

Suppose that for some $t_0\in [0,1]$ we have $\sigma_d(A(t_0))\not=\emptyset$, otherwise the theorem is obviously true.
Then it follows that there exist 
$$\Delta_j,\,\,\,\la_j(\cdot),\,\,\, E_j(\cdot)\,\,\,\text{ and } \,\,\,x_k^j(\cdot)$$
as in Lemma \ref{l:evf} such that
$\Delta_j\cap[0,1]\not=\emptyset$ for some $j\in\N$. 
By $\frakK$ denote the set of all $j$ such that $\la_j(t)\in (a,b)$ 
for some $t\in\Delta_j\cap [0,1]$ and for $j\in\frakK$ define
$$
\wt\Delta_j := \left\{t\in\Delta_j\cap [0,1] : \la_j(t)\in (a,b)\right\} = \la_j^{-1}((a,b))\cap [0,1].
$$
Due to \eqref{e:DGL} and the continuity of $\la_j(\cdot)$ the set $\wt\Delta_j$ is a (non-empty) 
subinterval of $\Delta_j$ which is open in $[0,1]$. For $j\in\frakK$, $t\in [0,1]$ and $k\in\{1, ..., m_j\}$ we set
\begin{align}
\begin{split}\label{e:tildelambda}
\wt{\la}_j(t) := 
\begin{cases} 
\lim_{s\downto\inf\wt\Delta_j} \la_j(s), \quad &0\le t\le\inf\wt\Delta_j,\\
\la_j(t), \quad &t\in\wt\Delta_j,\\
\lim_{s\upto\sup\wt\Delta_j} \la_j(s), \quad &\sup\wt\Delta_j\le t\le 1,
\end{cases}
\end{split}
\end{align}
\begin{align*}
\begin{split}
\wt{E}_j(t) := 
\begin{cases} 
E_j(t), \quad &t\in\wt\Delta_j,\\
0, \quad &t\in [0,1]\setminus\wt\Delta_j,
\end{cases}
\end{split}
\end{align*}
and
\begin{align*}
\wt{x}_k^j(t) := 
\begin{cases} 
x_k^j(t), \quad &t\in\wt\Delta_j,\\
0, \quad &t\in [0,1]\setminus\wt\Delta_j.
\end{cases}	
\end{align*}
The functions $\wt{\la}_j(\cdot)$, $\wt{E}_j(\cdot)$, and $\wt{x}_k^j(\cdot)$ are differentiable in all but at most two points 
$t\in [0,1]$ and for each $j\in\frakK$  the differential equation 
\begin{equation}\label{e:DGL2}
\wt\la_j'(t) = \frac{1}{m_j}\sum_{k=1}^{m_j}\left[C\wt x_k^j(t), \wt x_k^j(t)\right]\,\ge\,0
\end{equation}
holds in all but at most two points $t\in [0,1]$; cf.
\eqref{e:DGL}. In addition, the projections 
$\wt{E}_j(t)$ are $\product$-selfadjoint for every $t\in [0,1]$. The rest of this proof is divided into several steps.

\vskip 0.2cm\noindent
{\it 1. Basis representations:} By $E_C$ denote the spectral function of the non-ne\-ga\-tive operator $C$. Since $0$ is 
not a singular critical point of $C$, the spectral projections $E_C(\R^+)$, $E_C(\R^-)$ and $E_C(\{0\})$ exist. 
In particular, $E_C(\{0\})\calK = \ker C^2 = \ker C$ is a Krein space. Let
\begin{equation*}
\ker C = \calH_+\,[\ds]\,\calH_-
\end{equation*}
be an arbitrary fundamental decomposition of $\ker C$. Then with $\calK_\pm := \calH_\pm\,[\ds]\,E_C(\R^\pm)\calK$ 
we obtain a fundamental decomposition
$$
\calK = \calK_+\,[\ds]\,\calK_-
$$
of $\calK$. By $J$ denote the fundamental symmetry associated with this fundamental decomposition and set 
$\hproduct := [J\cdot,\cdot]$. Then $\hproduct$ is a Hilbert space scalar product on $\calK$, and $C$ is a 
selfadjoint operator in the Hilbert space $(\calK,\hproduct)$. By $\|\cdot\|$ denote the norm induced by $\hproduct$. 
Let $(\gamma_l)$ be an enumeration of the non-zero eigenvalues of $C$ (counting multiplicities). Since $C\in\frakS_p(\calK)$, we have
\begin{equation}\label{e:gamma_l}
(\gamma_l)\in\ell^p.
\end{equation}
Let $\{\vphi_l\}_l$ be an $\hproduct$-orthonormal basis of $\ol{\ran C}$ such that $\vphi_l$ is an eigenvector 
of $C$ corresponding to the eigenvalue $\gamma_l$. Then we have $|[\vphi_l,\vphi_i]| = \delta_{li}$. In the following we do
not distinguish the cases $\dim\ran C<\infty$ and $\dim\ran C=\infty$, that is, $l=1,\dots,m$ for some $m\in\N$ and $l\in\N$,
respectively. 

Consider the basis representation of $v\in\ol{\ran C}$ with respect to $\{\vphi_l\}_l$. There exist 
$\alpha_l\in\C$ such that $v = \sum_{l}\alpha_l\vphi_l$. Therefore
$$
[v,\vphi_k] = \sum_{l}\alpha_l[\vphi_l,\vphi_k] = \alpha_k [\vphi_k,\vphi_k]\quad\text{and}\quad
v = \sum_{l}\frac{[v,\vphi_l]}{[\vphi_l,\vphi_l]}\,\vphi_l.
$$
Consequently, for $x = u+v$, $u\in\ker C$, $v\in\ol{\ran C}$, we have $[x,\vphi_l]=[v,\vphi_l]$, $[Cx,x] = [Cx,v]$ and hence
\begin{align}
\begin{split}\label{e:[Cx,x]}
[Cx,x] 
&= \left[Cx,\sum_{l}\frac{[x,\vphi_l]}{[\vphi_l,\vphi_l]}\,\vphi_l\right]
= \sum_{l}[Cx,\vphi_l]\frac{[\vphi_l,x]}{[\vphi_l,\vphi_l]} = \sum_{l}[x,C\vphi_l]\frac{[\vphi_l,x]}{[\vphi_l,\vphi_l]}\\ 
&= \sum_{l}[x,\gamma_l\vphi_l]\frac{[\vphi_l,x]}{[\vphi_l,\vphi_l]} = \sum_{l}\frac{\gamma_l}{[\vphi_l,\vphi_l]}\big|[x,\vphi_l]\big|^2
= \sum_{l}|\gamma_l|\big|[x,\vphi_l]\big|^2,
\end{split}
\end{align}
where the non-negativity of $C$ was used in the last equality, cf. \eqref{guteformel}. 
Let $j\in\frakK$ be fixed, $t\in\wt\Delta_j$ and $x\in\calK$. Then
$$
E_j(t)x
= \sum_{k=1}^{m_j}[E_j(t)x,x_k^j(t)] x_k^j(t) 
= \sum_{k=1}^{m_j}[x,E_j(t)x_k^j(t)] x_k^j(t)
= \sum_{k=1}^{m_j}[x,x_k^j(t)] x_k^j(t). 
$$
If $t\in [0,1]\setminus\wt\Delta_j$, then $\wt{E}_j(t) = 0$ and $\wt{x}_k^j(t) = 0$, $k=1, ..., m_j$. Hence
\begin{equation}\label{e:kleine Summe}
\wt{E}_j(t)x = \sum_{k=1}^{m_j} [x,\wt{x}_k^j(t)]\wt{x}_k^j(t)
\end{equation}
holds for all $t\in [0,1]$ and all $x\in\calK$.

\vskip 0.2cm\noindent
{\it 2. Norm bounds: }In the following we prove that the projections $\wt{E}_j(t)$ are uniformly bounded 
in $j\in\frakK$ and $t\in [0,1]$. For $x\in\calK$ we have $\wt{E}_j(t)x\in E_{A(t)}([a,b])\calK$, and with Lemma \ref{l:lmm} we obtain
\begin{align*}
\|J\wt{E}_j(t)x\|\|x\|
&\ge (J\wt{E}_j(t)x,x) 
= [\wt{E}_j(t)x,x] 
= [\wt{E}_j(t)x,\wt{E}_j(t)x]\\
& \ge\delta\|\wt{E}_j(t)x\|^2
= \delta\|J\wt{E}_j(t)x\|^2.
\end{align*}
This implies
\begin{equation}\label{e:Absch JE_j}
\|J\wt{E}_j(t)\|\le\frac{1}{\delta}.
\end{equation}
Similarly, $\|E_{A(t)}((a,b))\|\le 1/\delta$ is shown to hold for $t\in [0,1]$. Consequently, the eigenvalues of $J\wt{E}_j(t)$ do not exceed $1/\delta$, and from $\dim J\wt{E}_j(t)\calK\le m_j$ it follows that the operator $J\wt{E}_j(t)$ has at most $m_j$ non-zero eigenvalues. Hence, its trace $\operatorname{tr}(J\wt{E}_j(t))$ satisfies
\begin{equation*}
\operatorname{tr}(J\wt{E}_j(t))\le\frac{m_j}{\delta}. 
\end{equation*}

\vskip 0.2cm\noindent
{\it 3. The main estimate: }Let $j\in\frakK$. For $t\in [0,1]$ we have
\begin{equation*}
\{\wt\la_j(t) : j\in\frakK,\,\wt\Delta_j\ni t\} = (a,b)\cap\sigma_d(A(t)) =: \Xi(t),
\end{equation*}
and it follows from the (strong) $\sigma$-additivity of the spectral function $E_{A(t)}$ (see, e.g., \cite{lmm}) that for every $x\in\calK$
\begin{equation}\label{e:Projektoridentitaet}
\sum_{j\in\frakK}\wt{E}_j(t) x \:
= \sum_{j\in\frakK,t\in\wt\Delta_j}\!\!\!\!E_j(t)x\;
= \sum_{\la\in\Xi(t)}E_{A(t)}(\{\la\})x
= E_{A(t)}\big((a,b)\big)x.
\end{equation}
From the differential equation \eqref{e:DGL2} we obtain for $j\in\frakK$
\begin{align}
\begin{split}\label{e:HDI}
\wt{\la}_j(1) - \wt{\la}_j(0) 
=\;\,\,\!& \frac{1}{m_j} \int_{0}^{1} \sum_{k=1}^{m_j}\left[C\wt{x}_k^j(t),\wt{x}_k^j(t)\right]\,dt\\
\stackrel{\eqref{e:[Cx,x]}}{=}&\frac{1}{m_j}\int_{0}^{1}\sum_{k=1}^{m_j}\sum_l|\gamma_l|\left|\left[\wt{x}_k^j(t),\vphi_l\right]\right|^2\,dt\\
=\;\,\,\!&\sum_l\:\frac{|\gamma_l|}{m_j}\int_{0}^{1}\left[\:\sum_{k=1}^{m_j}\left[\vphi_l,\wt{x}_k^j(t)\right]\wt{x}_k^j(t),\vphi_l\right]\,dt\\
\stackrel{\eqref{e:kleine Summe}}{=}&\sum_l\:\frac{|\gamma_l|}{m_j}\int_{0}^{1}\left[\wt{E}_j(t)\vphi_l,\vphi_l\right]\,dt.
\end{split}
\end{align}
For $j\in\frakK$ and $l$ we set
$$
\sigma_{jl} := \frac{1}{m_j}\int_{0}^{1}\left[\wt{E}_j(t)\vphi_l,\vphi_l\right]\,dt\quad \text{ and }\quad\sigma_j := \sum_l\sigma_{jl}.
$$
Then $\sigma_j\ge 0$ for all $j\in\frakK$, as $\sigma_{jl}\ge 0$ for all $l$. In fact, 
we have $\sigma_j > 0$ for each $j\in\frakK$. Indeed, if $\sigma_j = 0$ for some $j\in\frakK$, then for every $t\in [0,1]$
$$
\tr\big(J\wt{E}_j(t)\big)
= \sum_l\left(J\wt{E}_j(t)\vphi_l,\vphi_l\right) = \sum_l\left[\wt{E}_j(t)\vphi_l,\vphi_l\right] =  0,
$$
which implies $J\wt{E}_j(t) = 0$ (and thus $\wt{E}_j(t) = 0$), since the $\hproduct$-selfadjoint operator $J\wt{E}_j(t)$ 
has only non-negative eigenvalues. Therefore, $\wt\Delta_j = \emptyset$, which is not possible. Moreover,
\begin{align}
\begin{split}\label{e:erste Hilfsabsch}
\sigma_j  =& \frac{1}{m_j}\int_{0}^{1}\sum_l\left[\wt{E}_j(t)\vphi_l,\vphi_l\right]\,dt
					= \frac{1}{m_j}\int_{0}^{1}\sum_l\left(J\wt{E}_j(t)\vphi_l,\vphi_l\right)\,dt\\
					=& \frac{1}{m_j}\int_{0}^{1}\tr\big(J\wt{E}_j(t)\big)\,dt \;
					\le \;\frac{1}{m_j}\int_{0}^{1}\frac{m_j}{\delta}\,dt \;
					= \;\frac{1}{\delta}.
\end{split}
\end{align}
In addition (cf.\ \eqref{e:Absch JE_j} and \eqref{e:Projektoridentitaet}), for each $l$ we have
\begin{align}
\begin{split}\label{e:zweite Hilfsabsch}
\sum_{j\in\frakK}m_j\sigma_{jl}
&= \sum_{j\in\frakK}\int_{0}^{1}\left[\wt{E}_j(t)\vphi_l,\vphi_l\right]\,dt 
= \int_{0}^{1}\biggl[\,\sum_{j\in\frakK}\wt{E}_j(t)\vphi_l,\vphi_l\biggr]\,dt\\
&= \int_{0}^{1}\left[E_{A(t)}\big((a,b)\big)\vphi_l,\vphi_l\right]\,dt 
\le \int_{0}^{1}\bigl\|E_{A(t)}\bigl((a,b)\bigr)\bigr\|\|\vphi_l\|^2\,dt
\le \frac{1}{\delta}.
\end{split}
\end{align}
Let $j\in\frakK$. For $n\in\N$ we set $c_n := \sum_{l=1}^{n}\sigma_{jl}/\sigma_j\le 1$. Then the convexity
of $x\mapsto |x|^p$, \eqref{e:HDI}, and \eqref{e:erste Hilfsabsch} imply
\begin{align*}
\left|\wt{\la}_j(1) - \wt{\la}_j(0)\right|^p
&= \lim_{n \to \infty}c_n^p\left(\sum_{l=1}^{n}\frac{\sigma_{jl}}{c_n\sigma_j} \sigma_{j}|\gamma_l|\right)^p
\le \lim_{n \to \infty}c_n^{p-1}\sum_{l=1}^{n}\frac{\sigma_{jl}}{\sigma_j}\sigma_{j}^p|\gamma_l|^p\\
&\le \sum_{l=1}^\infty\sigma_{jl} \;\sigma_j^{p-1}|\gamma_l|^p 
\le \frac{1}{\delta^{p-1}}\sum_{l=1}^\infty\sigma_{jl}|\gamma_l|^p
\end{align*}
in the case that $\ran C$ is infinite dimensional (that is, $l=1,\dots\infty$); otherwise the above estimate 
holds with a finite sum on the right hand side. 
Hence, \eqref{e:zweite Hilfsabsch} and \eqref{e:gamma_l} yield
\begin{equation}\label{e:lp}
\sum_{j\in\frakK} m_j\left|\wt{\la}_j(1) - \wt{\la}_j(0)\right|^p
\le \frac{1}{\delta^{p-1}}\sum_{j\in\frakK}\sum_l m_j\sigma_{jl} |\gamma_l|^p
\le \frac{1}{\delta^{p}}\sum_l |\gamma_l|^p < \infty.
\end{equation}

\vskip 0.2cm\noindent
{\it 4. Final conclusion: } It suffices to consider the case $[a,b]\cap\sess(A) \not= \emptyset$, as otherwise $\sigma_p(A)\cap (a,b)$ and 
$\sigma_p(B)\cap (a,b)$ are finite sets and hence the theorem holds. We consider the following three possibilities separately: 
$a,b\in\sess(A)$, exactly one endpoint of $(a,b)$ belongs to $\sess(A)$, and $a,b\not\in\sess(A)$.

\vskip 0.1cm
(i) Assume that $a,b\in\sess(A)$. Then, by Lemma~\ref{l:evf} and \eqref{e:tildelambda} for all $j\in\frakK$ 
the values $\wt\la_j(0)$ and $\wt\la_j(1)$ either are boundary points of $\sess(A)=\sess(B)$ (see \eqref{sessz})
or points in the discrete spectrum of $A$ and $B$, respectively. Taking into account the multiplicities of the discrete
eigenvalues of $A$ and $B$ it is easy to construct 
sequences 
$$(\alpha_n)\subset\{\wt\la_j(0):j\in\frakK\}\quad\text{and}\quad(\beta_n)\subset\{\wt\la_j(1):j\in\frakK\}$$ 
such that $(\alpha_n)$ and  $(\beta_n)$ are extended enumerations of discrete eigenvalues of $A$ and $B$ in $(a,b)$
and $(\beta_n - \alpha_n)\in\ell^p$ by \eqref{e:lp}.

\vskip 0.1cm
(ii) Suppose that $a\notin\sess(A)$ and $b\in\sess(A)$ (the case $a\in\sess(A)$ and $b\notin\sess(A)$ is treated analogously). 
Then for each $j\in\frakK$ the value $\wt\la_j(1)$ is either a boundary point of $\sess(B)$ or a discrete eigenvalue of $B$. 
Hence, the sequence $(\beta_n)$ in (i) is an extended enumeration of discrete eigenvalues of $B$ in $(a,b)$. But it might 
happen that there exist indices $j\in\frakK$ such that $\wt\la_j(0) = a$, which is not a boundary point of 
$\sess(A)$ and not a discrete eigenvalue of $A$ in $(a,b)$. In the following we shall show that the number of such indices is finite. 
Then we just replace the corresponding values $\wt\la_j(0)$ in $(\alpha_n)$ by a point in $\partial\sess(A)\cap (a,b]$
and obtain an extended enumeration $(\alpha_n)$ of discrete eigenvalues of $A$ in $(a,b)$ such that $(\beta_n - \alpha_n)\in\ell^p$.

Assume that $\wt\la_j(0)=a$ for all $j$ from some infinite subset $\frakK_a$ of $\frakK$.
Then $\wt\la_j(t) = a$ for all $t\in [0,t_j]$, where $t_j := \inf\wt\Delta_j$, $j\in\frakK_a$. 
Observe that $a\in\sigma_d(A(t_j))$ (cf. Lemma~\ref{l:evf}) and $\la_j(t_j) = a$, and as $a\not\in\sess(A(t))$ for all $t\in[0,1]$, the set $\{t_j:j\in\frakK_a\}$
is an infinite subset of $[0,1]$. Hence we can assume that $t_j$ converges to some $t_0$, $t_j\not= t_0$ for all $j\in\frakK_a$,
and that the functions $\lambda_j$ are not constant.
Choose $\veps > 0$ such that $a - \veps > 0$ and
\begin{equation*}
\bigl([a-\veps,a)\cup(a,a+\veps]\bigr)\cap\sigma(A(t_0)) = \emptyset.
\end{equation*}
Either $t_0\not\in\Delta_j$ or $t_0\in\Delta_j$, in which case $|\la_{j}(t_0) - a|>\veps$ holds.
As $\lambda_j(t_j)=a$ for
each $j$ there exists $s_j$ between $t_0$ and $t_j$ such that $|\la_{j}(s_j) - a| = \veps$. Therefore, 
there exists $\xi_j$ between $s_j$ and $t_j$ such that
$$
\veps = |\la_j(t_j) - \la_j(s_j)| = \la_j'(\xi_j)|t_j - s_j|\le \la_j'(\xi_j)|t_j - t_0|.
$$
Hence, $\la_j'(\xi_j)\to\infty$ as $j\to\infty$. On the other hand, by
Lemma \ref{l:lmm} there exists $\delta_0 > 0$ such that $[x,x]\,\ge\,\delta_0\|x\|^2$ for all $x\in E_{A(t)}([a-\veps,\infty))\calK$
and $t\in [0,1]$. Together with \eqref{e:DGL} this implies
$$
\la_j'(\xi_j)\,\le\,\frac{\|C\|}{m_j}\sum_{l=1}^{m_j}\|x_j^l(\xi_j)\|^2\,\le\,\frac{\|C\|}{m_j\delta_0}
\sum_{l=1}^{m_j}[x_j^l(\xi_j),x_j^l(\xi_j)] = \frac{\|C\|}{\delta_0},
$$
a contradiction. Hence there exist at most finitely many $j\in\frakK$ such that $\wt\lambda_j(0)=a$.

\vskip 0.1cm
(iii) If $a,b\not\in\sess(A)$, we choose $c\in(a,b)\cap\sess(A)$ and construct the extended enumerations $(\alpha_n)$ and $(\beta_n)$
as the unions of the extended enumerations in $(a,c)$ and $(c,b)$, which exist by (ii).
\end{proof}



\section*{Contact information}
Jussi Behrndt: Institut f\"ur numerische Mathematik, Technische Universit\"at Graz, Steyrergasse 30, 8010 Graz, Austria, behrndt@tugraz.at

\vspace{0.4cm}\noindent
Leslie Leben: Institut f\"ur Mathematik, Technische Universit\"at Ilmenau, Postfach 10 05 65, 98684 Ilmenau, Germany, leslie.leben@tu-ilmenau.de

\vspace{0.4cm}\noindent
Friedrich Philipp: Institut f\"ur Mathematik, MA 6-4, Technische Universit\"at Berlin, Stra\ss e des 17.\ Juni 136, 10623 Berlin, Germany, fmphilipp@gmail.com
\end{document}